\documentstyle[12pt]{article}
\newcommand{\sect}[1]{\setcounter{equation}{0}
\section{#1}}
\newcommand{\subsect}[1]{\subsection{#1}}

\renewcommand{\theequation}{\arabic{section}.
\arabic{equation}}

\def\be{\begin{equation}}
\def\ee{\end{equation}}
\def\bea{\begin{eqnarray}}
\def\eea{\end{eqnarray}}

\def\p{\varepsilon}
\def\pp{\tilde P}
\def\kk{\tilde K}
\def\jj{{\tilde J}}
\def\hh{{\tilde H}}
\def\mm{{\tilde M}}

\def\>#1{{\bf #1}}
\def\1{\'{\i}}
\def\R{\rm I\kern-.2em R}
\def\back{\!\!\!\!\!\!}
\def\wpq{W_+^q}
\def\wcq{W_-^q}
\def\wp{W_+}
\def\wc{W_-}

\begin{document}
\thispagestyle{empty}
\hfill\today
\vspace{2cm}

\begin{center} {\LARGE{\bf{A New ``Null-Plane" Quantum
Poincar\'e Algebra}}} \end{center}

\bigskip\bigskip\bigskip

\begin{center} A. Ballesteros$^{1}$, F.J. Herranz$^{1}$,
M.A. del Olmo$^{2}$ and M. Santander$^{2}$ \end{center}

\begin{center} {\it { {}$^{1}$ Departamento de F\1sica
Aplicada, Universidad de Burgos} \\   E-09003, Burgos,
Spain}

{\it  { {}$^{2}$ Departamento de F\1sica Te\'orica,
Universidad de Valladolid}\\ E-47011, Valladolid, Spain}
\end{center}

\bigskip

\bigskip\bigskip\bigskip

\begin{abstract}

A new quantum deformation, which we call null-plane, of
the (3+1) Poin\-ca\-r\'e algebra is obtained. The
algebraic properties of the classical null-plane description are
generalized to this quantum deformation. In particular,
the classical isotopy subalgebra of the null-plane is
deformed into a Hopf subalgebra, and deformed spin
operators having classical commutation rules can be
defined. Quantum Hamiltonian, mass and position operators
are studied, and  the null-plane evolution is expressed in
terms of a deformed Schr\"odinger equation.

\end{abstract}

\newpage
%%%%%%%%%%%%%%%%% INTRODUCTION %%%%%%%%%%%
\sect {Introduction}

A central problem of a Hamiltonian formulation of quantum
relativistic systems consists in the obtention of a set of
dinamically independent variables in terms of the
generators of the (3+1) Poincar\'e group. For the
null-plane evolution \cite{Dirac} this has been thorougly
studied, and in this case it is well known that the
information provided by the Poincar\'e invariance into a
kinematical and a dynamical part. To achieve this by
preserving causality, one can consider that the initial
state of the system is defined on a light-like plane
$\Pi_n^\tau$ (a hypersurface of points $x$ of the
Minkowski space such that $n\cdot x=\tau$, where $n$ is a
light-like vector and $\tau$ a real constant). The
stability group $S_+$ of $\Pi_n^\tau$ gives the kinematics
inside the null-plane, and the dynamics is obtained by
means of the remaining generators (the Hamiltonians) that
map $\Pi_n^\tau$ into some other surface and describe the
evolution of the system. A brief summary of this
construction is given in Section 2 for the case $n=(\frac
12,0,0,\frac 12)$ (see \cite{LSAnn} for a general review).

In this letter we present a new quantum deformation of the
Poincar\'e algebra such that the classical null-plane
isotopy subalgebra $S_+$ is promoted into a Hopf
subalgebra under deformation. The classical null-plane
separation between kinematics and dynamics is therefore
preserved in the quantum case. So, the study of
deformations of the classical models built within this
framework (for instance, the infinite momentum frame
approach \cite{KSRep}, gauge field theory quantized on
null-planes \cite{KSop} and applications in Hadron
spectroscopy \cite{LSHad,HHLS}) could be physically
interesting.

We also recall that previous Hopf algebra deformations of
the (3+1) Poincar\'e algebra \cite{LRN,CD,BH4D} have been
obtained within the purely kinematical framework where
their classical counterparts live. However,
the quantum algebra here introduced induces a (dynamical)
deformation of the Schr\"odinger equation governing the
evolution of the wave functions with support on the
null-plane.

{}From an algebraic point of view the (3+1) Poincar\'e
deformation introduced in Section 3 is a quantization of a
coboundary Lie bialgebra generated by a (non-standard)
$r$--matrix that fulfills the classical Yang--Baxter
equation  (CYBE). Starting from a similar (2+1) null-plane
deformation, which was introduced in \cite{Bey}, this  new
(3+1) quantum Poincar\'e algebra is obtained by applying a
deformation embedding method (fully described in
\cite{BH4D,BHND} for the standard deformation).  It is
worth remarking that the deformation of the most relevant
structures in the classical null-plane description is
straightforward:  a centrally extended (2+1) Galilean
subalgebra of ${\cal P}(3+1)$ is transformed into a
quantum algebra under deformation, the quantum central
elements can be used to define the deformed mass and
Pauli--Lubanski operators, and besides the latter induces
a natural definition for the quantum spin operators, whose
commutation rules are shown to be the classical ones.
Quantum analogues of reduced Hamiltonians, position
operators and the deformed uncertainty relations linked to
them are studied in Section 4. Some relevant remarks close
the paper.

%%%%%%%%%%% SECTION 2 %%%%%%%%%%%%%%%%%%%%%%%
\sect {Classical Poincar\'e algebra in a null-plane basis}

Hereafter we shall consider the null-plane ``orthogonal"
to the light-like vector $n=(\frac 12,0,0,\frac 12)$ as
our initial surface \cite{LSAnn}. A coordinate  system
well adapted to null-planes  $\Pi_n^\tau$ is given by
\be
x^- = n\cdot x =\frac 12(x^0-x^3)=\tau,\quad
x^+=x^0+x^3,\quad x_T=(x^1,x^2). \label{aa}
\ee
Hence, a point $x$ contained in the null-plane $\Pi_n^\tau$
will be labelled by the coordinates $(x^+, x^1,x^2)$. The
remaining coordinate $x^-$ plays the role of a time. A
particular null-plane is $\Pi_n^0, (\tau=0)$, i.e., $x^- =
n\cdot x =0$; this plane is mapped onto itself by the
boosts generated by $K_3$, which leave the transverse
coordinates $x_T=(x^1,x^2)$ unchanged and $e^{\chi K_3}$
maps $x_+$ on $e^{\chi}x_+$. The Poincar\'e algebra
generators can be classified according to their
transformation property with respect to $K_3$; any
generator $A$ obeying
\be
 [K_3,A]=\gamma A,\label{ab}
\ee
is referred to as an operator of ``goodness" $\gamma$. The
ten generators of  ${\cal P}(3+1)$ (adapted to the
null-plane coordinates in terms of the usual physical
basis) can be listed in relation to their goodness in the
following way:
\bea
\gamma=+1:&&\quad E_1=\frac 12(K_1+J_2), \quad
E_2=\frac 12(K_2 - J_1),\quad P_+=\frac 12 (P_0+P_3);\cr
\gamma=0:&&\quad\qquad  K_3,\qquad J_3,\qquad P_1,\qquad
P_2;\label{ac}\\ \gamma=-1:&&\quad   F_1=\frac 12(K_1
-J_2), \quad F_2=\frac 12(K_2 + J_1),\quad P_-=\frac 12
(P_0-P_3).\nonumber
\eea

The non-vanishing Lie brackets of ${\cal P}(3+1)$ are
$(i,j=1,2)$:
\bea
 && [K_3,P_+]=P_+,\quad
[K_3,P_-]=-P_-,\quad [K_3,E_i]=E_i,\quad
[K_3,F_i]=-F_i,\cr
&& [J_3,P_i]=-\p_{ij3}P_j, \qquad  [J_3,E_i]=-\p_{ij3}E_j,
\qquad [J_3,F_i]=-\p_{ij3}F_j,\cr &&
[E_i,P_j]=\delta_{ij}P_+,\qquad  [F_i,P_j]=\delta_{ij}P_-,
\qquad [P_+,F_i]=-P_i,\label{com2}\\ &&
[E_i,F_j]=\delta_{ij}K_3 +\p_{ij3}J_3,\qquad
[P_-,E_i]=-P_i.\nonumber
\eea

The set of all Poincar\'e generators with the same
goodness $\gamma$ spans a subgroup $G_\gamma$ (in
particular, the subgroups $G_{\pm 1}$ are Abelian).
Moreover, there exist two seven-parametrical Poincar\'e
subgrups with a semidirect product structure:
$S_\pm=G_0\odot G_{\pm 1}$.  The stability group of the
plane $\Pi_n^0$ coincides with the group $S_+$. The
remaining three generators, which originate $G_{-1}$, act
on $\Pi_n^0$ as follows: the generator $P_-$ translates
$\Pi_n^0 $ into $\Pi_n^\tau$, while the generators $F_i$
rotate it around the surface of the light-cone $x^2=0$.
Therefore, if $x^-=\tau$ is considered as an evolution
parameter, then $P_-$ and $F_i$ describe the dynamical
evolution  from the null-plane $x^-=0$; they are called
Hamiltonians.

The two Casimirs of ${\cal P}(3+1)$ are given by:
\bea
&& M^2=2 P_-P_+- P_1^2 - P_2^2,\label{cas1}\\ &&
W^2=W_{13}^2+W_{23}^2+ \wp\wc +\wc\wp , \label{cas2}
\eea
where
\bea
&& W_{13}=K_3 P_1 +E_1 P_- -F_1 P_+,\quad
W_{23}=K_3 P_2 +E_2 P_- -F_2 P_+, \cr
&& \wp=E_1 P_2 -E_2 P_1 +J_3 P_+ , \qquad
\wc=F_1 P_2 -F_2 P_1 +J_3 P_- .\label{comp}
\eea
The second-order Casimir $M^2$ is the square of the mass
operator ($M$) and the fourth-order one $W^2$ is the
square of the  Pauli-Lubanski operator. Lie brackets
involving the components $W_i$ and the null-plane
generators of  ${\cal P}(3+1)$ are given in Appendix A.

%%%%%%%%%%%%%% SECTION 3 %%%%%%%%%%%%%%%%%%%%%%%%%

\sect {Quantum Poincar\'e algebra in a null-plane basis}

In \cite{BH4D}, the standard deformation of
${\cal P}(3+1)$ has been derived by using the so-called
{\it deformation embedding method}, a kind of quantum
analogue of the classical embedding of the ${\cal P}(2+1)$
algebra into ${\cal P}(3+1)$ (see also Ref. \cite{BHND} for
arbitrary dimension). In this way, the $(2+1)$ dimensional
deformation of a given (affine) algebra turns out to be
the cornerstone for the obtention of its higher
dimensional generalizations. Since a null-plane quantum
${\cal P}(2+1)$ algebra has been given in \cite{Bey}, the
same procedure can be now applied as follows:

1) Let us look for classical subalgebras of
${\cal P}(3+1)$ isomorphic to the null-plane algebra
${\cal P}(2+1)$: there exist two of them. Namely,
\be
\Pi_{13}=\langle K_3,E_1,F_1,P_1,P_+,P_-\rangle ,
\qquad  \Pi_{23}=\langle
K_3,E_2,F_2,P_2,P_+,P_-\rangle.\label{pib}
\ee

2) We assume that $U_z\Pi_{13}$ and
$U_z\Pi_{23}$, endowed with the null-plane $(2+1)$
deformation of \cite{Bey} are ``restrictions" \cite{BHND}
of the quantization for $U_z{\cal P}(3+1)$.

3) Consider the simplest Ansatz for the
coproduct and commutation rules of $U_z{\cal P}(3+1)$
consistent with these restrictions.

4) Check the selfconsistency of the structure
obtained in 3) by imposing that the coproduct is a Hopf
algebra homomorphism.

The final result is the following null-plane deformation of
$U_z{\cal P}(3+1)$, with coproduct given by:
\bea
&&\Delta(X)=1\otimes X+X\otimes 1, \qquad \mbox{for}\quad
X\in\{P_+,E_1,E_2,J_3\},\cr
&&\Delta(Y)=e^{-zP_+}\otimes Y+Y\otimes e^{zP_+},\quad
\mbox{for}\quad Y\in\{P_-,P_1,P_2\},\cr
&& \Delta(F_1)=e^{-zP_+}\otimes F_1
+F_1\otimes e^{zP_+} +ze^{-zP_+} E_1 \otimes P_- -
zP_-\otimes E_1 e^{zP_+}\cr
 &&\qquad\qquad  +ze^{-zP_+} J_3 \otimes P_2 -
z P_2\otimes J_3 e^{zP_+},\label{co3}\\
&& \Delta(F_2)=e^{-zP_+}\otimes F_2+F_2\otimes e^{zP_+}
+ze^{-zP_+} E_2 \otimes P_- - zP_-\otimes E_2 e^{zP_+}\cr
&&\qquad\qquad  +ze^{-zP_+} J_3 \otimes P_1 -
z P_1\otimes J_3 e^{zP_+},\cr
&& \Delta(K_3)=e^{-zP_+}\otimes K_3+K_3\otimes e^{zP_+}
+ze^{-zP_+} E_1 \otimes P_1 - zP_1\otimes E_1 e^{zP_+}\cr
&&\qquad\qquad +ze^{-zP_+} E_2 \otimes P_2 -
z P_2\otimes E_2 e^{zP_+}; \nonumber
\eea
commutation relations:
\bea
&& [K_3,P_+]=\frac {\sinh zP_+}z,\quad [K_3,P_-]=-P_-
\cosh zP_+,\quad [K_3,E_i]=E_i\cosh zP_+,\cr
&& [K_3,F_1]=- F_1 \cosh zP_+  +  z   E_1 P_-\sinh zP_+
- z^2P_2\,\wpq,\cr && [K_3,F_2]=- F_2 \cosh zP_+
+  z  E_2 P_- \sinh zP_+  + z^2P_1\,\wpq,\cr
&& [J_3,P_i]=-\p_{ij3}P_j, \qquad [J_3,E_i]=-\p_{ij3}E_j,
\qquad [J_3,F_i]=-\p_{ij3}F_j,\cr
&& [E_i,P_j]=\delta_{ij}\frac {\sinh zP_+}z, \qquad
[F_i,P_j]=\delta_{ij}P_-\cosh zP_+,\label{brack3}\\
&& [E_i,F_j]=\delta_{ij}K_3 +\p_{ij3}J_3\cosh zP_+,\qquad
[P_+,F_i]=-P_i,\cr && [F_1,F_2]=z^2 P_- \wpq +
z  P_- J_3 {\sinh zP_+}, \qquad [P_-,E_i]=-P_i.\nonumber
\eea
(the quantum component $\wpq$ of the Pauli--Lubanski vector
is defined below (\ref{qcas5}));
counit and  antipode:
\be
\epsilon(X) =0;\quad \gamma(X)=-e^{ {3z} P_+}\ X\
e^{- {3z}P_+},\quad   \mbox{for $X\in
\{P_\pm,P_i,E_i,F_i,J_3\}$}. \label{anti3}
\ee

As a byproduct of the embedding method, the underlying
coboundary Lie bialgebra is also obtained. Explicitly, the
classical $r$--matrix
\be
r=2z(K_3 \wedge P_+ +E_1\wedge P_1
+ E_2\wedge P_2),\label{rclas}
\ee
is a skew-solution of the CYBE (it can be identified with
one appearing in Zakrzewski's classification
\cite{Zakr}), and generates the cocommutator
$\delta (X)=[1\otimes X + X\otimes 1, r]$
that provides the first order term in $z$ within the
coproduct (\ref{co3}):
\bea
&& \delta(X)=0,\qquad
\mbox{for}\quad X\in\{P_+,E_1,E_2,J_3\},\cr
&&\delta(Y)=2z Y\wedge P_+ ,\quad  \mbox{for}\quad
Y\in\{P_-,P_1,P_2\},\cr &&\delta(F_1)=2z(F_1\wedge P_+
+E_1\wedge P_- +   J_3\wedge P_2),\label{coco3}\\
&&\delta(F_2)=2z(F_2\wedge P_+ +E_2\wedge P_- +
J_3\wedge P_1),\cr
&&\delta(K_3)=2z(K_3\wedge P_+ +E_1\wedge P_1 +
E_2\wedge P_2).\nonumber
\eea

The mass (\ref{cas1}) and the Pauli-Lubanski (\ref{cas2})
operators can be deformed as follows:
\be
M^2_q=2 P_-\frac{\sinh zP_+}z - P_1^2 -
P_2^2,\label{qcas3}
\ee
\be
W^2_q=(W_{13}^q)^2+(W_{23}^q)^2 +
\cosh zP_+\left\{\wpq\wcq +\wcq\wpq\right\} -
{z^2} M_q^2 (\wpq)^2, \label{qcas4}
\ee
where
\bea
&&
W_{13}^q= K_3 P_1 +E_1 P_-\cosh zP_+
-F_1\frac{\sinh zP_+}z,\cr &&W_{23}^q= K_3 P_2
+E_2 P_-\cosh zP_+ -F_2\frac{\sinh zP_+}z,\cr
&& \wpq=  E_1 P_2 - E_2 P_1 +J_3\frac {\sinh zP_+}z
,\label{qcas5}\\
&&\wcq=  F_1 P_2 - F_2 P_1 +J_3 P_-{\cosh zP_+}. \nonumber
\eea
The commutation rules between the components $W_i^q$ and
the Hopf algebra generators are given in Appendix A; they
show that $W^2_q$ belongs to the center of
$U_z{\cal P}(3+1)$ as $M^2_q$.

It is interesting to analyse now the deformed algebras
given by the generators of the classical subgroups
$G_\gamma,\ \gamma=-1,0,+1$ and $S_\pm$. It is a rather
remarkable fact that the generators of the null-plane
stability group, $S_+$, determine a deformed Hopf
subalgebra of $U_z{\cal P}(3+1)$ (that we shall denote by
$U_z S_+$), while the Lie subalgebra generating $G_{+1}$
gives a trivial (i.e., non deformed) Hopf subalgebra under
deformation. On the other hand, $G_0$, $G_{-1}$ and $S_-$
do not originate Hopf subalgebras. Moreover, although the
classical group $G_{-1}$ is Abelian,  this is not the case
in the quantum version.

There is another physically interesting deformed Hopf
subalgebra in $U_z{\cal P}(3+1)$: the set of generators
$\langle J_3,E_1,E_2,P_1,P_2,P_+,P_-\rangle$  generates a
quantum deformation, $U_z\widetilde{\cal G}(2+1)$, of the
(2+1) extended Galilean algebra. The generator $P_+$ is
central, and the structure of $U_z\widetilde{\cal G}(2+1)$
can be studied by considering a new physical basis
$\{\hh, \pp_1,\pp_2,\kk_1,\kk_2,\jj\}$ linked, as usual, to
the time translation, space translations, boosts and space
rotation, respectively.  For the Galilean case we
introduce a central element $\mm$ which can be identified
with the mass of a free particle. Explicitly, the
relations between the null-plane and the usual kinematical
bases for this subalgebra are:
\be
\back\back\back U_z\widetilde{\cal G}(2+1):\
\langle J_3,E_1,E_2,P_1,P_2,P_-,P_+\rangle\equiv
\langle -\jj_3,\kk_1,\kk_2, \pp_1,\pp_2,\hh,\mm\rangle.
\label{bn}
\ee

The Hopf structure of $U_z\widetilde{\cal G}(2+1)$ has the
following coproduct and non-vanishing  commutation rules:
\bea
&&\Delta(X)=1\otimes X+X\otimes 1,\quad
X\in\{\mm,\jj_3,\kk_1,\kk_2\},\cr &&\Delta(Y)=
e^{-  z \mm}\otimes Y + Y\otimes e^{  z \mm} ,\quad
Y\in\{\hh,\pp_1,\pp_2\}. \label{bo}
\eea
\bea
&&[\jj_3,\kk_i]=\p_{ij3} \kk_j,\qquad
[\jj_3,\pp_i]=\p_{ij3} \pp_j,\cr
&&[\kk_i,\hh]=\pp_i,\qquad
[\kk_i,\pp_j]=\delta_{ij}\frac{\sinh z\mm}z .\label{bp}
\eea
The center of $U_z\widetilde{\cal G}(2+1)$ is
generated by the mass $\mm$, the energy $E_q$ and the
intrinsic angular momentum $L_q$:
\bea
&&E_q^2=\pp_1^2+\pp_2^2  - \hh \frac {\sinh   z  \mm }z
,\label{bq}\\ &&L_q= \kk_1\pp_2-\kk_2\pp_1 -
\jj_3 \frac{\sinh   z  \mm } z  .\label{br}
\eea

%%%%%%%%%%%%%%%%%%%%%%% SECTION 4 %%%%%%%%%%%%%%%%

\sect {Spin, Hamiltonians and Position Operators}

In the classical Poincar\'e algebra the spin
$\vec{\cal S}$ satisfies the commutation rules:
\be
[{\cal S}_i,{\cal S}_j]=\p_{ijk} {\cal S}_k ,\qquad
[M,{\cal S}_i]=0,\quad i,j,k=1,2,3,\label{da}
\ee
where the components ${\cal S}_i$ are given in terms of the
Pauli--Lubanski vector by
\be
M{\cal S}_1=W_{23} +(P_1/P_+)\wp,\quad M{\cal S}_2=
W_{13} -(P_2/P_+)\wp,\quad {\cal S}_3=\wp / P_+ .\label{db}
\ee
The spin $\vec{\cal S}$ commutes with all the
generators of the stability group $S_+$ except $J_3$. In
particular, the component ${\cal S}_3$ is a central
element of $S_+$. Its eigenvalue $h$ is called the
null-plane helicity.

Relations (\ref{db}) enable us to write  the Hamiltonians
$P_-$, $F_1$ and $F_2$ in terms of the mass and spin
operators:
\bea
&& P_-=[M^2+P_1^2+P_2^2](1/2P_+),\cr
&& F_1=[K_3P_1+E_1P_- -(M{\cal S}_2 +
P_2 {\cal S}_3)](1/P_+),\label{dc}\\ && F_2=[K_3P_2+E_2P_-
-(M{\cal S}_1 - P_1 {\cal S}_3)](1/P_+).\nonumber
\eea
In these expressions the most relevant objects
(rather than ${\cal S}_1$, ${\cal S}_2$ and $M$) are the
products $M{\cal S}_1$, $M{\cal S}_2$ and $M^2$. These
quadratic elements are called ``reduced Hamiltonians".

\subsect {Quantum Spin and Hamiltonians}

The quantum counterpart of this classical structure can be
straightforwardly obtained once we know the
$q$--Casimirs of  $U_z{\cal P}(3+1)$
(\ref{qcas3}--\ref{qcas5}). We define the components
${\cal S}_i^q$ of the quantum spin in such a way that the
following relations hold:
\be
M_q{\cal S}^q_1=W_{23}^q +\frac{zP_1}
{\tanh zP_+}\wpq,\quad  M_q{\cal S}^q_2=W_{13}^q
-\frac{zP_2}{\tanh zP_+}\wpq,\quad {\cal S}^q_3=
\frac{z\wpq}{\sinh zP_+} .\label{dd}
\ee
Relations (\ref{aaa}--\ref{aae}) of Appendix A show that
the  $q$--spin  $\vec{\cal S}^q$, so introduced,  verifies
the following commutation rules:
\be
[{\cal S}_i^q,{\cal S}_j^q]=\p_{ijk} {\cal S}_k^q ,\qquad
[M_q,{\cal S}_i^q]=0.\label{de}
\ee
Similarly to the classical case, $\vec{\cal S}^q$ commutes
with all the generators of the quantum stability group
$U_zS_+$ except  $J_3$; furthermore, the component
${\cal S}_3^q$ (quantum helicity) is a central element for
$U_zS_+$.

The Hamiltonians now adopt the following deformed
expressions:
\bea
&& P_-=[M^2_q+P_1^2+P_2^2]\frac{z}{2\sinh zP_+} ,\cr
&& F_1=[K_3P_1+E_1P_-\cosh zP_+ -(M_q{\cal S}_2^q +
P_2\cosh zP_+ {\cal S}_3^q)]\frac{z}
{ \sinh zP_+},\label{df}\\ && F_2=[K_3P_2+E_2P_-\cosh zP_+
-(M_q{\cal S}_1^q - P_1 \cosh zP_+{\cal S}_3^q)]
\frac{z}{\sinh zP_+},\nonumber
\eea
and we call the products $M_q{\cal S}_1^q$,
$M_q{\cal S}_2^q$ and $M^2_q$ quantum reduced Hamiltonians.
Note that in the quantum case $[F_1,F_2]\neq 0$.

A quantum differential realization for $U_z{\cal P}(3+1)$
can be obtained in the momentum representation with
coordinates $\>p=(p_+,p_1,p_2)=(p_+,p_T)$ and $p_+ >0$.
The simplest procedure to derive it consists in getting
firstly a realization for the quantum stability  group
$U_zS_+$. It is easy to check that it reads
\bea
&& P_+=p_+,\qquad P_i=p_i,\qquad E_i=\frac{\sinh zp_+}{z}
{\partial_i},\cr
&& J_3= p_1   {\partial_2} - p_2 {\partial_1}+
{\cal S}_3^q,\qquad  K_3=\frac{\sinh zp_+}z  {\partial_+},
\ \  i=1,2, \label{dg}
\eea
where $\partial_i=\frac{\partial}{\partial p_i}$.  By
introducing (\ref{dg}) in the expressions of the quantum
Hamiltonians (\ref{df}) we get
\bea
&& P_-=\frac{z( M^2_q+p_T^2)}{2\sinh zp_+} ,\cr
&& F_1= p_1 {\partial_+}  +\frac {z(M_q^2+p_T^2)}
{2\tanh zp_+}\partial_1- \frac z{\sinh zp_+}
(M_q{\cal S}^q_2 + p_2\cosh zp_+{\cal S}_3^q),\label{dh}\\
&& F_2= p_2 {\partial_+}   +\frac {z(M_q^2+p_T^2)}
{2\tanh zp_+}\partial_2- \frac z{\sinh zp_+}
(M_q{\cal S}^q_1 - p_1\cosh zp_+{\cal S}_3^q),\nonumber
\eea
with $p_T^2=p_1^2+p_2^2$.

Let us consider a basis vectors $|\>p,h\rangle$ for this
representation, i.e.,
$P_l|\>p,h\rangle=p_l|\>p,h\rangle,\ l=+,1,2$, and
$S_3|\>p,h\rangle=h|\>p,h\rangle$. The action of the
remaining operators is given by (\ref{dg}) and (\ref{dh}).
The norm of these eigenvectors is
\be
\langle p^\prime, h|p,h\rangle= (2\pi)^3 2
\frac{\sinh zp_+}{z}\delta (p^\prime -p)
\delta _{h^\prime,h}. \label{dha}
\ee
Now, we introduce a quantum inner product of two arbitrary
states depending on the three momentum coordinates
$\>p=(p_+,p_1,p_2)$ and the helicity $h$ (the quantum spin
has classical commutation rules):
\be
\langle \phi|\psi\rangle =\frac 1{(2\pi)^3}\sum_{h}
\int \frac{z\, d^3\>p}{2\sinh zp_+} \phi_h^\ast(\>p)
\psi_h(\>p),\label{inner}
\ee
where $\phi_h(\>p)= <p,h|\phi>$. Note that in the limit
$z \to 0$ we recover the appropiate classical inner
product. According to expression  (\ref{inner}), generators
$P_+$, and $P_l$ are hermitian operators while $K_3$,
$J_3$ and $E_l$ are skew-hermitian. Therefore,
\be
\hat P_+ \equiv P_+,\quad \hat P_l \equiv P_l,\quad
\hat K_3 \equiv i K_3,\quad \hat E_l \equiv i E_l, \quad
\hat J_3 \equiv i J_3,\quad (l=1,2),\label{herm}
\ee
gives a hermitian representation of $U_zS_+$ relative to
(\ref{inner}) which induces a hermitian representation of
whole $q$--algebra $U_z({\cal P}(3+1)$. Note that the
commutations rules for the hermitian operators
(\ref{herm}) are also given by (3.3) by adding the
imaginary unit $i$ in the r.h.s of the corresponding
commutators. Only $P_-$ out of the remaining three
generators is hermitian. On the other hand, let us also
remark that in order to extend the $*$--operation at the
Hopf algebra level $z$ has to be consider a dimensional
real parameter $[z]=[p_+]^{-1}$ (see coproduct (3.2)).

We can consider the coordinate $x^-$ as a ``time"
($\tau$), and the complete wave function will be
$\tau$--dependent $f_h(\>p,\tau)$. The Hamiltonian
$\hat P_-\equiv P_-$ provides the evolution of the
null-plane in terms of a $q$--deformed Schr\"odinger
equation $i{\partial_\tau}f_h(\>p,\tau)=\hat
P_-f_h(\>p,\tau)$. From the representation (\ref{dh}) we
obtain
\be
i{\partial_\tau}\psi_h(\>p,\tau)=
\frac{z( m^2_q+p_T^2)}{2\sinh zp_+} \psi_h(\>p, \tau),
\label{ddd}
\ee
where the functions $\psi_h(\>p, \tau)$
are square integrable under (\ref{inner})
\be
\frac 1{(2\pi)^3}\sum_{h}\int \frac{z\, d^3\>p}
{2\sinh zp_+} |\psi_h(\>p, \tau)|^2 < +\infty.
\label{os}
\ee
The r.h.s. of this $q$--Schr\"odinger equation can be seen
as a deformation of the kinetic term
$H_0=(m^2_q+p_T^2)/2p_+$ of the null-plane bound state
equation in quantum chromodynamics \cite{KSop,HHLS}. In
fact, when $z$ is small enough, we can write
\be
\hat P_- =\frac{z( m^2_q+p_T^2)}{2\sinh zp_+} =
\frac{m^2_q+p_T^2}{2 p_+} - \frac{1}{12}\, z^2 \,p_+
\,(m^2_q+p_T^2) + o(z^4) = H_0 + V_z, \label{qcd}
\ee
and the null-plane deformation intrinsically includes a
dynamical part $V_z$ to be computed together with the
quantum chromodynamics interaction.

\subsect {Quantum Position Operators}

In the classical case there exist two null-plane (or mean
transverse) position operators $\hat Q_i$   $(i=1,2)$ (see
\cite{LSAnn}) which measure the position of the center of
the longitudinal momentum. Explicitly,
\be
\hat Q_i=  \hat  E_i/\hat P_+,\qquad i=1,2;\label{di}
\ee
and they satisfy the following classical commutation rules
with the generators of the stability group $S_+$:
\bea
&& [\hat Q_i,\hat E_j]=0,\qquad [\hat Q_i,
\hat K_3]=0,\qquad  [\hat Q_i,\hat J_3]=i\p_{ij3}
\hat  Q_j,\cr
&& [\hat Q_i,\hat P_+]=0,\qquad [\hat Q_i,\hat P_j]= i
\delta_{ij}, \qquad [\hat Q_1,\hat Q_2]=0. \label{dj}
\eea
{}From the definition of $\hat Q_i$ (\ref {di}) in the
classical version of the representation (\ref {dg}) the
position operators are
$\hat Q_i=i\frac{\partial}{\partial p_i}$. They commute
with the reduced Hamiltonians $M^2$, $M{\cal S}_i$ and
with the helicity ${\cal S}_3$. The transverse velocity
operator ($\frac{d}{d \tau}\hat Q_i$) is given by the
commutator of $\hat Q_i$ with the Hamiltonian $\hat
P_-$:    \be  \frac{d}{d\tau}\hat Q_i=[\hat Q_i,\hat
P_-]=i  \hat  P_i/\hat P_+ .\label{dk}
\ee

In order to find quantum deformations for the position
operators, $\hat Q_i$,  we make the following Ansatz
\be
\hat Q_i= \hat E_i/f(z,\hat P_+), \qquad i=1,2,\label{dl}
\ee
where $f(z,\hat P_+)$ is an arbitrary (smooth)
function such that
 $\lim_{z\to 0}f(z,\hat P_+)=\hat P_+$ and
$f(z,p_+)\neq 0$ for $z$ real and $p_+>0$. Commutators
(\ref{dj}) and (\ref{dk}) become in the quantum case (we
omit the arguments in the function $f$):
\bea
&& \back\back\back [\hat Q_i,\hat E_j]=0,\quad
[\hat Q_i,\hat K_3]=i\hat E_i\left\{\frac {f'}{f^2}
\frac{\sinh z\hat P_+}{z} -
\frac {\cosh z\hat P_+}{f}\right\},\quad
[\hat Q_i,\hat J_3]=i \p_{ij3} \hat Q_j,\cr
&& \back\back\back  [\hat Q_i,\hat P_+]=0,\quad
[\hat Q_i,\hat P_j]= i\delta_{ij}
\frac{\sinh z\hat P_+}z\frac 1{f},\quad
[\hat Q_1,\hat Q_2]=0, \quad [\hat Q_i,\hat P_-]=
i\frac{\hat P_i}f,  \label{dm}
\eea
where $f'$ is ``the formal derivative" of $f$ with respect
to $\hat P_+$. Moreover, whatever the choice of the
function $f$ be, the operators $\hat Q_i$ are hermitian
with respect to (\ref{inner}), and they commute with the
quantum reduced Hamiltonians and with the quantum helicity
--similarly to the classical case. We analyze the most
natural possibilities:

\noindent (1) $f=\sinh z\hat P_+/z$. The only deformed
bracket in (\ref{dm}) is the one corresponding to the
transverse velocity operator, i.e.,
$[\hat Q_i,\hat P_-]= i z{\hat P_i}/\sinh z\hat P_+$. Thus,
the quantum realization for $\hat Q_i$ is
\be
\hat Q_i=\frac{z\hat E_i}{\sinh z\hat P_+},\label{dn}
\ee
which, in terms of the differential realization (\ref{dg})
provides the undeformed result
$\hat Q_i \equiv i\frac{\partial}{\partial p_i}$ for the
position operators. As a result, the ``generalized"
uncertainty principle derived from (\ref{dn}) is given by
\be
\Delta \hat Q_i\Delta \hat P_i \geq \frac12,
\label{dmm}
\ee
(in this case $[\hat Q_i,\hat P_j]=i\delta_{ij}$, since we
have taken $\hbar=c=1$).

\noindent   (2) $f=\tanh z\hat P_+/z$. Now the commutation
rules (\ref{dm}) read:
\bea
&& \back\back\back  [\hat Q_i,\hat E_j]=0,\quad
[\hat Q_i,\hat K_3]=i\hat E_i\, z\sinh z\hat P_+, \quad
[\hat Q_i,\hat J_3]=i\p_{ij3}\hat  Q_j, \cr
&& \back\back\back     \quad [\hat Q_i,\hat P_+]=0,\quad
[\hat Q_i,\hat P_j]=i\delta_{ij}  \cosh z\hat P_+,\quad
[\hat Q_1,\hat Q_2]=0.\label{do}
\eea
The time derivative of $\hat Q_i$ turns into
\be
\frac{d}{d\tau}\hat Q_i= [\hat Q_i,\hat P_-]=
i \frac{z\hat P_i}{\tanh z\hat P_+}, \label{dnb}
\ee
and the corresponding generalized uncertainty relation is
given by
\be
\Delta \hat Q_i\Delta \hat P_i \geq
\frac12\Delta(\cosh z\hat P_+). \label{dnc}
\ee
A power series expansion gives
\be \Delta \hat Q_i\Delta \hat P_i \geq
\frac12(1+\frac{z^2}{2}(\Delta \hat P_+)^2+\dots),
\label{dnd}
\ee
a result similar to the one presented in \cite{Mag} that
concerns a generalized uncertainty principle in Quantum
Gravity. In our case the ``classical" value
$\Delta \hat Q_i\Delta \hat P_i \geq \frac 12$ appears
deformed by addition of terms depending on
$\Delta \hat P_+$.

%%%%%%%%%%%% CONCLUDING REMARKS %%%%%%%%%%%%
\sect {Concluding remarks}

We have presented a new quantum deformation of the algebra
${\cal P}(3+1)$ whose natural framework is the null-plane
basis for this algebra. From the point of view of quantum
group theory, this deformation presents some interesting
properties that distinguish it from the other ones already
known. The essential point is that the classical
$r$--matrix (\ref{rclas}) is triangular; thus, a
$\ast_\hbar$--product quantizing the algebra of
representative functions on the group exists. First order
information concerning the corresponding null-plane
quantum Poincar\'e group can be extracted from the Lie
bialgebra structure (\ref{coco3}). In particular, for the
(now non-commutative) coordinates
$(\hat p_1,\hat p_2,\hat p_+)$ we shall have
\be
[\hat p_1,\hat p_+]=2\,z\,\hat p_1 + o(z^2),\qquad
[\hat p_2,\hat p_+]=2\,z\,\hat p_2 + o(z^2).
\ee
Moreover, since $r$ is a solution of the CYBE, a universal
$R$--matrix turning this non-standard deformation into a
triangular Hopf algebra could be hopefully derived. All
these formal aspects deserve further study.

{}From a physical point of view, the most relevant feature
of the results here displayed is the parallelism between
the classical subalgebras and the Hopf algebra ones
induced after deformation. This fact is the cornerstone
behind the straightforward generalization of the classical
constructions  intrinsically containing a modification of
both the kinematical and dynamical classical symmetries.
In this way, by means of a hermitian representation of the
quantum isotopy algebra of the null-plane $\Pi_n^0$, it has
been possible to define a family of (hermitian)
Newton--Wigner position operators corresponding to the
transverse coordinates $q_1$ and $q_2$. These operators
are not univocally defined as there exists some freedom in
the choice of a function $f(z,\hat P_+)$  given in the
above section. The Lie bracket between the position and
the momentum operators gives rise to a $q$--deformed
generalized uncertainty principle whose expression is
\be
\Delta \hat Q_i\Delta \hat P_i \geq \frac 12\Delta
\left (\frac{\sinh z\hat P_+}{zf(z,\hat P_+)}\right ).
\ee
Note that expressions (\ref{dmm}) and (\ref{dnc}) are
particular cases of (5.2) according to the different
choices of $f$.

The time evolution of a null-plane deformed system (given
by the $q$--Schr\"odinger equation (\ref{ddd})) is also a
consecuence of the existence of an extended (2+1) Galilean
Hopf subalgebra, in which $\hat P_-$ plays the role of the
deformed Hamiltonian and $\hat P_+$ is the mass operator
in perfect agreement with the classical version. Quantum
deformed reformulations of classical models based in the
null-plane description (and in particular, the role of the
deformed coproduct for composite systems) seem to be worth
studying.

\newpage

%%%%%%%%%%%%%%%% APPENDIX %%%%%%%%%%%%%%%%%%%%%

\noindent {\large{\bf{Appendix A: Commutation relations
for the classical and quantum Pauli--Lubanski operator}}}

\appendix

\setcounter{equation}{0}

\renewcommand{\theequation}{A.\arabic{equation}}

\bigskip

The non-vanishing commutation relations of the components
$W_i$ of the Pauli-Lubanski vector (\ref{comp}) among
themselves and  with the generators of ${\cal P}(3+1)$ are:
\bea
&&\back [J_3,W_{13} ]=-W_{23},\quad [E_2,W_{13}]=\wp,\quad
[F_2,W_{13} ]=-\wc ,\cr
&&\back [J_3,W_{23}]=W_{13},\qquad [E_1,W_{23} ]=-\wp,
\quad [F_1,W_{23}]=\wc,\\
&&\back [E_1,\wc ]=W_{23},\quad [E_2,\wc ]=-W_{13} ,\quad
[K_3,\wc ]= -\wc,\cr
&&\back [F_1,\wp ]=-W_{23},\qquad [F_2,\wp ]=W_{13} ,\quad
[K_3,\wp ]=\wp;\cr && \cr
&&\back [W_{13},\wp]=\wp P_1 + W_{23} P_+,\quad
[W_{13},\wc]=-\wc P_1 + W_{23} P_-,\cr
&&\back [W_{23},\wp]=\wp P_2 - W_{13} P_+, \quad
[W_{23},\wc]=-\wc P_2 - W_{13} P_-,\\
&&\back [W_{23},W_{13}]=\wp P_- + \wc P_+ , \quad
[\wp,\wc]= W_{13} P_1 + W_{23} P_2   .\nonumber
\eea

The quantum analogues of the above relations, which allow
to  prove that $W^2_q$ (\ref{qcas4}) belongs to the center
of  $U_z{\cal P}(3+1)$, are:
\bea
&&[J_3,W_{13}^q]=-W_{23}^q,\qquad  [F_1,W_{13}^q]= {z^2}
\wpq P_1P_2,\cr
&&[E_2,W_{13}^q]=\wpq\cosh zP_+,\qquad [K_3,W_{13}^q]=z
\wpq \, P_2\sinh zP_+ ,\label{aaa}\\
&&[F_2,W_{13}^q]=-\wcq\cosh zP_+ + {z^2} \wpq
\left\{ M_q^2 + P_2^2\right\},\cr
&&\cr
&&[J_3,W_{23}^q]=W_{13}^q,\qquad  [F_2,W_{23}^q]=
-\ {z^2} \wpq P_1P_2,\cr
&&[E_1,W_{23}^q]=-\wpq\cosh zP_+,\qquad [K_3,W_{23}^q]=
- z \wpq\, P_1\sinh zP_+ ,\\
&&[F_1,W_{23}^q]=\wcq\cosh zP_+ - {z^2} \wpq
\left\{M_q^2 + P_1^2\right\},\cr && \cr
&& [E_1,\wcq]=W_{23}^q,\qquad [F_1,\wcq]=- {z^2}
 \wpq P_1P_- ,\cr && [E_{2},\wcq]=- W_{13}^q,\qquad
[F_2,\wcq]=-{z^2} \wpq P_2P_- ,\\
&&[K_3,\wcq]=-\wcq\cosh zP_+ + {z^2} \wpq
\left\{M_q^2 - P_-\frac{\sinh zP_+}z\right\},\nonumber
\eea
\be
[F_1,\wpq]=- W_{23}^q,\quad  [F_{2},\wpq]=
W_{13}^q,\quad  [K_3,\wpq]= \wpq \cosh zP_+ ;
\ee
\bea
&&\back [W_{13}^q,\wpq]=\wpq P_1 \cosh zP_+
+ W_{23}^q \frac{\sinh zP_+}z,\cr
&&\back [W_{23}^q,\wpq]=\wpq P_2 \cosh zP_+
- W_{13}^q \frac{\sinh zP_+}z,\cr
&&\back [W_{23}^q,W_{13}^q]=\wpq \left\{ P_-\cosh^2 zP_+
- zM_q^2 \sinh zP_+ \right\} + \wcq
\frac{\sinh 2zP_+}{2z},\label{aae}\\
 &&\back [W_{13}^q,\wcq]=-\wcq P_1 \cosh zP_+ +
W_{23}^q P_-{\cosh zP_+}  +z^2 M_q^2 \wpq,\cr
&&\back [W_{23}^q,\wcq]=-\wcq P_2 \cosh zP_+ -
W_{13}^q P_-{\cosh zP_+}  +z^2 M_q^2 \wpq,\cr
&&\back  [\wpq,\wcq]= W_{13}^q P_1 + W_{23}^q P_2.
\nonumber
\eea

\newpage

\bigskip \medskip

\noindent {\large{{\bf Acknowledgements}}}

\bigskip This work has been partially supported by DGICYT
de Espa\~na (Projects   PB91--0196 and PB92--0255).

\end{document}